\documentclass[10pt]{article}
\usepackage{epsfig,amssymb,amsmath}
\usepackage{palatino}
\usepackage{verbatim}
\usepackage{enumerate}

\def\baa {\begin{eqnarray*}}
\def\eaa {\end{eqnarray*}}

\def \la {\lambda}
\def \al {\alpha}

\def \n {{\lfloor\frac{n}{2} \rfloor}}
\def \nn {{\lfloor\frac{n+1}{2} \rfloor}}

\def\Frac#1#2{\mbox{\large${\textstyle \frac{#1}{#2}}$}}

\def\Frac#1#2{\mbox{\large${\textstyle \frac{#1}{#2}}$}}

\def \ra {{\quad\Rightarrow\quad}}
\def \lr {{\quad\Leftrightarrow\quad}}

\def \tr {{\rm tr\,}}

\def\R{{\mathbb R}}
\def\N{{\mathbb N}}
\def\A{{\mathbf A}}
\def\B{{\mathbf B}}
\def\PP{{\cal P}}

\def\OO{{\cal O}}

\def\wt{\widetilde}
\def\wh{\widehat}

\def\la{\lambda}

\def \l {\ell}
\newtheorem{lemma}{Lemma}[section]
\newtheorem{proposition}[lemma]{Proposition}
\newtheorem{corollary}[lemma]{Corollary}
\newtheorem{theorem}[lemma]{Theorem}
\newtheorem{remark}[lemma]{Remark}

\newtheorem{definition}[lemma]{Definition}

\def\be  {\begin{equation}}
\def\ee  {\end{equation}}
\def\ba  {\begin{eqnarray}}
\def\ea  {\end{eqnarray}}
\def\baa {\begin{eqnarray*}}
\def\eaa {\end{eqnarray*}}
\def\bc  {}
\makeatletter
\@addtoreset{equation}{section}
\makeatother
\def\proof{\medskip\noindent{\bf Proof.} }
\def\qed{\hfill $\Box$}
\newcommand {\lb} {\label}
\newcommand {\rf[1]} {(\ref{#1})}
\newcommand {\ds} {\displaystyle}

\begin{document}

\title{On the Markov inequality in the $L_2$-norm with the Gegenbauer weight}

\author{G.\,Nikolov, A.\,Shadrin}

\date{}
\maketitle

\begin{abstract}
Let $w_\la(t) := (1-t^2)^{\la-1/2}$, where $\la > -\frac{1}{2}$, be
the Gegenbauer weight function, let $\|\cdot\|_{w_\la}$ be the
associated $L_2$-norm,
$$
    \|f\|_{w_\la} = \left\{\int_{-1}^1 |f(x)|^2 w_\la(x)\,dx\right\}^{1/2}\,,
$$
and denote by $\PP_n$ the space of algebraic polynomials of degree
$\le n$.
We study the best constant $c_n(\la)$ in the Markov inequality
in this norm
$$
   \|p_n'\|_{w_\la} \le c_n(\la) \|p_n\|_{w_\la}\,,\qquad p_n \in \PP_n\,,
$$
namely the constant
$$
   c_n(\la) := \sup_{p_n \in \PP_n} \frac{\|p_n'\|_{w_\la}}{\|p_n\|_{w_\la}}\,.
$$
We derive explicit lower and upper bounds for the Markov constant
$c_n(\la)$, which are valid for all $n$ and $\la$.
\end{abstract}

\smallskip


\section{Introduction}


Let $w_\la(t) := (1-t^2)^{\la-1/2}$, where $\la > -\frac{1}{2}$, be
the Gegenbauer weight function, let $\|\cdot\|_{w_\la}$ be the
associated $L_2$-norm,
$$
    \|f\|_{w_\la} = \left\{\int_{-1}^1 |f(x)|^2 w_\la(x)\,dx\right\}^{1/2}\,,
$$
and denote by $\PP_n$ the space of algebraic polynomials of degree
$\le n$.
In this paper, we study the best constant $c_n(\la)$ in the Markov inequality
in this norm
\begin{equation}\label{e1.1}
   \|p_n'\|_{w_\la} \le c_n(\la) \|p_n\|_{w_\la}\,,\qquad p_n \in \PP_n\,,
\end{equation}
namely the constant
$$
   c_n(\la) := \sup_{p_n \in \PP_n} \frac{\|p_n'\|_{w_\la}}{\|p_n\|_{w_\la}}\,.
$$
Our goal is to derive {\it good} and {\it explicit}
lower and upper bounds for the Markov constant $c_n(\la)$
which are valid for {\it all} $n$ and $\la$, i.e., to find constants
$\underline{c}(n,\la)$ and $\overline{c}(n,\la)$ such that
$$
    \underline{c}(n,\la) \le c_n(\la) \le \overline{c}(n,\la)\,,
$$
with a small ratio $\frac{\overline{c}(n,\la)}{\underline{c}(n,\la)}$.

It is known that, for a fixed $\la$, $c_n(\la)$ grows like
$\OO(n^2)$, and that the asymptotic value
$$
  c_*(\la) := \lim_{n\to\infty} \frac {c_n(\la)}{n^2}
$$
is equal to $1/(2\,j_{\frac{2\la-3}{4}})$, with $j_{\nu}$ being the
first positive zero of the Bessel function $J_{\nu}$, see
\cite[Thms. 1.1--1.3]{BD2010}, whereby it can be shown that
$c_*(\la)$ behaves like $\OO(\la^{-1})$. There is also a number of
more precise results.

For $\la = \frac{1}{2}$ (the constant weight $w_{\frac{1}{2}} \equiv
1$), it follows from the Schmidt result \cite{s44} that
$$
    \frac{1}{\pi} (n+\Frac{3}{2})^2
\le c_n(\Frac{1}{2})
\le \frac{1}{\pi}(n+2)^2\,.
$$
For $\la = 0,1$ (the Chebyshev weights $w_0(x) =
\frac{1}{\sqrt{1-x^2}}$ and $w_1(x) = \sqrt{1-x^2}$, respectively),
Nikolov \cite{n03} proved that \be \lb{n}
\begin{array}{l}
  0.472135 n^2 \le c_n(0) \le 0.478849 (n+2)^2\,, \\[0.5ex]
  0.248549 n^2 \le c_n(1) \le 0.256861 (n+\Frac{5}{2})^2\,.
\end{array}
\ee In \cite{ans16}, we obtained an upper bound valid for all $n$
and $\la$, \be \lb{l_1}
   c_n(\la) \le \frac{(n+1)(n+2\la +1)}{2\sqrt{2\la+1}}\,,
\ee however, the already mentioned asymptotics $c_*(\la) =
\OO(\la^{-1})$ shows that this result is not optimal.

The main result of this paper is lower and upper bounds for $c_n(\la)$ which
are uniform with respect to $n$ and $\la$. They show, in particular, that
$$
   [c_n(\la)]^2 \asymp \frac{1}{\la^2}n(n+2\la)^3\,.
$$
For $n=1,2$ the exact values of the Markov constant are easily
computable:
\be \lb{e0}
   [c_1(\la)]^2 = 2(1+\la), \qquad
   [c_2(\la)]^2 = \frac{4(2+\la)(2+2\la)}{2\la+1}\,.
\ee Therefore, we consider below the case $n \ge 3$. Our main result
is

\begin{theorem} \lb{thm1}
For all $\la > -\frac{1}{2}$ and $n\ge 3$, the best constant $c_n(\la)$
in the Markov inequality
$$
   \|p_n'\|_{w_\la} \le c_n(\la) \|p_n\|_{w_\la}\,,\qquad p_n \in \PP_n\,,
$$
admits the estimates
\ba
   \frac{1}{4}\frac{n^2(n+\la)^2}{(\la+1)(\la+2)} \;
 < & [c_n(\la)]^2 &
 < \;\frac{n(n+2\la+2)^3}{(\la+2)(\la+3)},
       \qquad\quad \la \ge 2\,; \lb{e1} \\
   \frac{(n+\la)^2(n+2\la')^2}{(2\la+1)(2\la+5)} \;
 < & [c_n(\la)]^2 &
 < \;\frac{(n+\la+\la''+2)^4}{2(2\la+1)\sqrt{2\la+5}},
       \qquad \la > -\Frac{1}{2}\,, \lb{e2}
\ea where $\la' = \min\,\{0,\la\}$, $\,\la'' = \max\,\{0,\la\}$.
\end{theorem}

As a consequence, we can specify the following bounds for the
asymptotic value $c_*(\la)$:

\begin{corollary}
For any $\la > -\frac{1}{2}$, the asymptotic Markov constant
$c_*(\la) = \lim\limits_{n\to\infty}n^{-2} c_n(\la)$ satisfies the
inequalities
$$
\frac{1}{(2\la+1)(2\la+5)} < [c_*(\la)]^2 <
 \begin{cases}\,
\ds{\frac{1}{2(2\la+1)\sqrt{2\la+5}}}\,,&\quad -\frac{1}{2}<\la\leq
 \la^{*}\,,\vspace*{1mm}\\ \, \ds{\frac{1}{(\la+2)(\la+3)}}\,, &\quad \la>\la^{*}\,,
 \end{cases}
$$
where $\la^{*}\approx 25$\,.
\end{corollary}

The lower bound in \rf[e1] follows from that in \rf[e2] and is less
accurate, we put it in this form to make the comparison between the
two bounds in \rf[e1] more obvious.

The upper bound in \rf[e2] does not have the right order
$\OO(n^4/\la^2)$ in $\,\la\,$ (for $\la$ fixed), however this bound
serves not only for the case $-\frac{1}{2} <\la < 2$, but for a
fixed $\la\in [2,\la^{*}]$ and $n\geq n_0(\la)$ it is also better
than the one in \rf[e1].

In the next corollary, we set $\la = 0,1$ in the upper estimate
\rf[e2], and that improves the upper estimates in \rf[n] for the
Chebyshev weights. When coupled with the lower estimate from \rf[n],
this gives rather tight bounds.

\begin{corollary}
For the Chebyshev weights $w_0(x) = \frac{1}{\sqrt{1-x^2}}$ and
$w_1(x) = \sqrt{1-x^2}$, we have
\baa
  0.472135\, n^2 \le c_n(0) \le 0.472871\, (n+2)^2\,, &&\\
  0.248549\, n^2 \le c_n(1) \le 0.250987\, (n+4)^2\,. &&
\eaa
\end{corollary}

The lower and upper estimates in \rf[e1] have different
orders with respect to $\la$. However we can get a perfect match
with slightly less accurate constants.

\begin{theorem} \lb{thm5}
For all $\la \ge 7$ and $n \ge 3$, the best constant $c_n(\la)$ in
the Markov inequality satisfies \be \lb{e4}
   \frac{1}{16} \frac{n(n+2\la)^3}{\la^2}
\le [c_n(\la)]^2
\le \frac{n(n+2\la)^3}{\la^2}\,.
\ee
\end{theorem}

\begin{corollary} \lb{thm6}
For the Markov constant $c_n(\la)$ we have the following asymptotic
estimates:

i)~~ $\displaystyle{ \sqrt{n}\;
 \le \lim_{\la\to\infty}\frac{c_n(\la)}{\sqrt{2\la}}
 \le  \sqrt{3n}}\,$;\bigskip

ii)~~ $\displaystyle (n -\Frac{1}{2})(n-1) \le  \lim_{\la\to
-\frac{1}{2}} c_n(\la) \cdot 2\sqrt{2\la + 1}
 \le  (n + \Frac{3}{2})^2$\,.
\end{corollary}

Part ii) follows from \rf[e2]. Though part i) does not formally
follow from Theorem \ref{thm5}, it follows from a part of its proof.

Let us describe briefly how these results are obtained.

It is well-known that the squared best constant in the Markov
inequality in the $L_2$-norm with arbitrary (and possibly different)
weights for $p$ and $p'$ is equal to the largest eigenvalue of a
certain positive definite matrix, in our case we have \be \lb{B}
    [c_n(\la)]^2 = \mu_{\max}(\B_n)\,,
\ee where the matrix $\B_n$ is specified in Sect. 2. We obtain then
lower and upper bounds for $\mu_{\max}(\B_n)$ using three values
associated with the matrix $\B_n$ and its eigenvalues $(\mu_i)$
(note that $\mu_i>0$):\smallskip

a) the trace
$$
   \tr(\B_n) := \sum b_{ii} = \sum \mu_i\,;
$$

b) the max-norm
$$
   \|\B_n\|_\infty = \max_i  \sum_j |b_{ij}|\,;
$$

c) the Frobenius norm
$$
   \|\B_n\|_F^2 := \sum_{i,j} |b_{ij}|^2 = \tr(\B_n \B_n^T) = \sum
   \mu_i^2\,.
$$

Clearly, we have
\be \lb{mu<}
   {\rm i)}  \quad \mu_{\max} \le \tr(\B_n)\,, \qquad
   {\rm ii)}\quad \mu_{\max} \le \|\B_n\|_\infty\,, \qquad
   {\rm iii)} \quad \mu_{\max} \le \|\B_n\|_F\,,
\ee and generally $\mu_{\max} \le \|\B_n\|_*$, where $\|\cdot\|_*$
is any matrix norm. The upper estimate \rf[l_1] cited from
\cite{ans16} is exactly the first inequality $\mu_{\max} \le
\tr(\B_n)$, and as we noted, this estimate is not optimal. The
better upper bounds \rf[e1]-\rf[e2] in Theorem \ref{thm1} are
obtained from (\ref{mu<}.ii) and (\ref{mu<}.iii), respectively.

For the lower bounds we use the inequalities
\be \lb{mu>}
    {\rm i')} \quad\mu_{\max}
\ge \frac{\sum \mu_i^2}{\sum \mu_i}
 =  \frac{\|\B_n\|_F^2}{\tr(\B_n)}\,, \qquad
    {\rm ii')} \quad \mu_{\max}(\B_n)
\ge \max_i b_{ii}\,. \ee Inequality (i') gives the lower estimates
in \rf[e1]-\rf[e2], and combination of (i') and (ii') yields the
lower bound in \rf[e4].

The paper is organised as follows. In Sect.\,\ref{pre}, following
our previous studies \cite{ans16}, we give an explicit form of the
matrix $\B_n$ appearing in \rf[B]. Sects.\,\ref{pre}-\ref{lemmas}
contain some auxiliary inequalities.
In Sect.\,\ref{infty}, we find an upper bound for the max-norm
$\|\B_n\|_\infty$, and in Sect\,.\ref{Fro} we give both lower and
upper estimates for the Frobinuis norm $\|\B_n\|_F$. Finally, in
Sect.\,\ref{main} we prove the upper and the lower estimates in
Theorems \ref{thm1}-\ref{thm5} using inequalities \rf[mu<]-\rf[mu>]
and relation \rf[B]. Here we have used the expression for
$\tr(\B_n)$ and for diagonal elements $b_{ii}$ found in
\cite{ans16}.

The formulas for the trace, the max-norm and the Frobenius norm of
a matrix are straightforward once the matrix elements are known,
so the main technical issues are, firstly, in finding reasonable
upper and lower bounds for the entries of the matrix $\B_n = (b_{ij})$
which are expressed initially in terms of the Gamma function $\Gamma$,
and, secondly, in finding reasonable estimates for their sums.
The first issue is dealt with in Sect.\,\ref{est}, where we show that
$$
       b_{jk}
\asymp \frac{f_\sigma(j)}{f_\tau(k)}\,, \qquad
       f_\al(x) = x^{\al_1}(x+\Frac{\la}{2})^{\al_2}(x+\la)^{\al_3}\,,
$$
and the second one in Sect.\,\ref{lemmas}\,, where we give elementary
but effective upper and lower bounds for the integrals of the type
$$
     \int_{x_0}^x f(t)\,dt, \qquad
f(x) = (x+\gamma_1)^{\al_1}(x+ \gamma_2)^{\al_2}
     \cdots(x+\gamma_r)^{\al_r}.
$$


\section{Preliminaries} \lb{pre}

In this section, we quote a result obtained earlier in \cite{ans16},
which equate the Markov constant $c_n(\la)$ with the largest eigenvalue
of a specific matrix $\B_n$.

\begin{definition} \rm
For $n \in \N$, set $m := \lfloor \frac{n+1}{2} \rfloor$ and define
symmetric positive definite matrices $\A_m,\wt\A_m \in \R^{m \times
m}$ with entries $a_{kj}$ and $\wt{a}_{kj}$ given by
\begin{equation}\lb{ab0}
     a_{kj}
:= \Big(\sum_{i=1}^{\min(k,j)} \al_i^2 \Big) \beta_k \beta_j\,, \qquad
     \wt a_{kj}
:= \Big(\sum_{i=1}^{\min(k,j)} \wt\al_i^2 \Big) \wt\beta_k \wt\beta_j\,,
\end{equation}
so that
\begin{equation}\lb{A_m}
   \A_m
:= \begin{pmatrix}
  \alpha_1^2\beta_1^2 &
  \alpha_1^2\beta_1\beta_2 & \alpha_1^2\beta_1\beta_3 & \cdots &
  \alpha_1^2\beta_1\beta_m \\
  \alpha_1^2\beta_1\beta_2 &
  \Big(\sum_{i=1}^{2}\alpha_i^2\Big)\beta_2^2 &
  \Big(\sum_{i=1}^{2}\alpha_i^2\Big)\beta_2\beta_3 & \cdots &
  \Big(\sum_{i=1}^{2}\alpha_i^2\Big)\beta_2\beta_m \\
  \alpha_1^2\beta_1\beta_3 &
  \Big(\sum_{i=1}^{2}\alpha_i^2\Big)\beta_2\beta_3 &
  \Big(\sum_{i=1}^{3}\alpha_i^2\Big)\beta_3^2 & \cdots &
  \Big(\sum_{i=1}^{3}\alpha_i^2\Big)\beta_3\beta_m \\
  \vdots & \vdots &
  \vdots & \ddots & \vdots \\
  \alpha_1^2\beta_1\beta_m &
  \Big(\sum_{i=1}^{2}\alpha_i^2\Big)\beta_2\beta_m &
  \Big(\sum_{i=1}^{3}\alpha_i^3\Big)\beta_3\beta_m & \cdots &
  \Big(\sum_{i=1}^{m}\alpha_i^2\Big)\beta_m^2
  \end{pmatrix}\,,
\end{equation}
with the same outlook for $\wt\A_m$. The numbers $\al_k, \beta_k$
and $\wt\al_k, \wt\beta_k$ are given by
\begin{eqnarray}
   \al_k := (2k-1 + \la) h_{2k-1},
& \beta_k := \Frac{1}{h_{2k}}\,; \lb{b}\\
   \wt\al_k := (2k-2 + \la) h_{2k-2},
& \wt\beta_k := \Frac{1}{h_{2k-1}}\,,\lb{wb}
\end{eqnarray}
where
\begin{equation}\lb{h_i}
   h_i^2
:= h_{i,\lambda}^2
:= \frac{\Gamma(i+2\lambda)}{(i+\lambda)\Gamma(i+1)}\,.
\end{equation}
Note that
\begin{equation}\lb{wbb}
   \wt\al_k = \al_{k-\frac{1}{2}}\,, \qquad
   \wt\beta_k = \beta_{k-\frac{1}{2}}\,.
\end{equation}
\end{definition}

\begin{definition} \rm
For $n \in \N$, set
\be \lb{AB0}
   \B_n := \left\{\begin{array}{ll}
   4 \A_m, & n = 2m; \\
   4 \wt\A_m, & n = 2m-1.
   \end{array} \right.
\ee
\end{definition}

\begin{theorem}[\cite{ans16}, Theorem 3.2]
Let $c_n(\la)$ be the best constant in the Markov inequality
\rf[e1.1]. Then
$$
    [c_n(\la)]^2 = \mu_{\max}(\B_n)\,,
$$
where $\mu_{\max}(\B_n)$ is the largest eigenvalue of the matrix $\B_n$.
\end{theorem}

\begin{remark} \rm
Appearance of two matrices $\A_m$ and $\wt\A_m$ reflects the fact
that the extreme polynomial $\wh p_n$ for the Markov inequality with
an even weight function $w(x) = w(-x)$ is either odd or even. The
latter is a relatively simple conclusion, what is not obvious though
is whether $\wh p_n$ is of degree exactly $n$ and not $n-1$. In
\cite{ans16}, we proved that for the Gegenbauer weights $w_\la$,
$$
   \mu_{\max}(\wt\A_m) <  \mu_{\max}(\A_m) <  \mu_{\max}(\wt\A_{m+1})
$$
and this implies that ${\rm deg}\,\wh p_n = n$, hence $[c_n(\la)]^2$
is the largest eigenvalue of $\A_m$ or $\wt\A_m$ for $n=2m$ or
$n=2m-1$, respectively.
\end{remark}

We finish this section by simplifying the expressions for $a_{kj}$
and thus for the matrix $\A_m$ as follows. From \rf[ab0], we derive
$$
     a_{kj}
:= \Big(\sum_{i=1}^{\min(k,j)} \al_i^2 \Big) \beta_k \beta_j
 = \left\{ \begin{array}{ll}
   \frac{\beta_k}{\beta_j}\,
      \big(\beta_j^2\sum_{i=1}^{j}\alpha_i^2\big), & j < k\,, \\
      \frac{\beta_j}{\beta_{k}}\,
      \big(\beta_k^2\sum_{i=1}^{k}\alpha_i^2\big), & j > k\,,
   \end{array}\right.
$$
so that
\be \lb{akj}
   a_{jj} = \beta_j^2\sum_{i=1}^{j}\alpha_i^2\,,\qquad
   a_{kj}
 = \left\{ \begin{array}{ll}
    \frac{\beta_k}{\beta_j}\,a_{jj}, & j < k\,, \\
    \frac{\beta_j}{\beta_{k}}\,a_{kk}, & j > k\,.
   \end{array}\right.
\ee  Respectively,
$$
   \A_m
= \begin{pmatrix}
a_{11} & \frac{\beta_2}{\beta_1} a_{11} & \frac{\beta_3}{\beta_1} a_{11} &
\cdots & \frac{\beta_m}{\beta_1} a_{11} \\
\frac{\beta_2}{\beta_1} a_{11} & a_{22} & \frac{\beta_3}{\beta_2} a_{22} &
\cdots & \frac{\beta_m}{\beta_2} a_{22} \\
\frac{\beta_3}{\beta_1} a_{11} & \frac{\beta_3}{\beta_2} a_{22} & a_{33} &
\cdots & \frac{\beta_m}{\beta_3} a_{33} \\
  \vdots & \vdots &
  \vdots & \ddots & \vdots \\
\frac{\beta_m}{\beta_1} a_{11} & \frac{\beta_m}{\beta_2} a_{22} &
\frac{\beta_m}{\beta_3} a_{33} &
\cdots & a_{mm}
\end{pmatrix}\,.
$$
Note that $\A_m$ and $\A_{m+1}$ are embedded. An analogous
representation and embedding hold for $\wt{\A}_m$.

\section{Estimates for $a_{kk}$ and $\frac{\beta_k}{\beta_j}$}  \lb{est}

We will need upper and lower estimates for the elements of matrices
$\A_m$ and $\wt\A_m$, namely
$$
   a_{kk} = \beta_k^2\sum_{i=1}^{k}\alpha_i^2\,,\qquad
   a_{kj}
 = \left\{ \begin{array}{ll}
    \frac{\beta_k}{\beta_j}\,a_{jj}, & j < k\,, \\
    \frac{\beta_j}{\beta_{k}}\,a_{kk}, & j > k\,.
   \end{array}\right.
$$
We found expression for $a_{kk}$ and $\widetilde{a}_{kk}$ in
\cite[Lemmas 2.1(ii) and 2.2(ii)]{ans16}, those are quoted in
Proposition \ref{sumab}, and in this section we obtain inequalities
for the ratios $\frac{\beta_k}{\beta_j}$.

\begin{proposition}[\cite{ans16}]\lb{sumab}
The following identities hold: \ba
    (i)\quad a_{kk} := \beta_k^2\,\sum_{i=1}^{k}\al_i^2
&=&  c_0 f_0(k), \lb{al} \\
    (ii)\quad \wt a_{kk} := \wt\beta_k^2\,\sum_{i=1}^{k}\wt\al_i^2
&=&  c_0 f_0(k-\Frac{1}{2}),  \lb{wtal} \ea where
$$
   c_0 := \frac{4}{2\la+1}, \qquad f_0(x) := x(x+\Frac{\la}{2})(x+\la)\,.
$$
\end{proposition}

\begin{proposition}\lb{bj/bk}
Let $j,k\in \N$, $j < k$. Then the coefficients $\beta_k$ in \rf[b]
satisfy the following relations:

\medskip
(i)\quad If $\,-\frac{1}{2} <\la\leq 0\,$ or $\,\la\ge 1$, then
\be \lb{b<}
    \Big(\frac{j}{k}\Big)^{2\la-2}
\le \frac{\beta_{k}^2}{\beta_j^2}
\le \Big(\frac{j+\la}{k+\la}\Big)^{2\la-2}\,.
\ee

(ii)\quad If  $\,0<\la\le 1$,  then
\be \lb{b>}
    \Big(\frac{j}{k}\Big)^{2\la-2}
\ge \frac{\beta_{k}^2}{\beta_j^2}
\ge \Big(\frac{j+\la}{k+\la}\Big)^{2\la-2}\,.
\ee
\end{proposition}

\proof Denote the left-hand, the middle and the right-hand side
terms in \rf[b<]-\rf[b>] by $\l(\la)$, $m(\la)$ and $r(\la)$,
respectively. From definitions \rf[b] and \rf[h_i] we have
\be
\lb{b1}
   m(\la) := \frac{\beta_{k}^2}{\beta_j^2}
 = \frac{\Gamma(2j+2\lambda)}{(2j+\lambda)\Gamma(2j+1)}
   \Big(\frac{\Gamma(2k+2\lambda)}{(2k+\lambda)\Gamma(2k+1)}\Big)^{-1}\,,
\ee
and using the functional equation $\Gamma(t+1) = t\,\Gamma(t)$ we see that
\be \lb{lmr}
   m(\la)
= \left\{\begin{array}{cl}
  \big(\frac{k}{j}\big)^{2}, & \la = 0, \\
  1, & \la = 1,
  \end{array} \right.
\ra \l(\la) = m(\la) = r(\la), \quad \la=0,1\,. \ee We shall prove
inequalities \rf[b<]-\rf[b>] for the logarithms of the values
involved.

1) Let us start with the proof of the left-hand side
inequalities in \rf[b<]-\rf[b>].
Consider the difference of the logarithms of the middle
and the left-hand side terms,
$$
   g(\la)
:= \log m(\la) - \log\l(\la)
 = \log m(\la) - (2\la-2)\log\frac{j}{k}
$$
We need to prove that $g(\la) \le 0$ for $\la \in [0,1]$ and that
$g(\la) > 0$ otherwise. Since $g(0)=g(1)=0$ by \rf[lmr], it suffices
to show that $g''(\la) > 0$ for all $\la > - \frac{1}{2}$, i.e.,
that $[\log m(\la)]'' > 0$.

From \rf[b1], we have
$$
   \log m(\la)
 = \log\Gamma(2j+2\lambda)-\log\Gamma(2k+2\lambda)
    -\log\frac{2j+\lambda}{2k+\lambda}
    -\log\frac{\Gamma(2j+1)}{\Gamma(2k+1)}\,,
$$
therefore, using the digamma function $\psi(t) :=
\Gamma'(t)/\Gamma(t)$, we obtain
$$
    [\log m(\la)]'
 = 2 \,\big[\psi(2j+2\la)-\psi(2k+2\la)\big]
     - \Big[\frac{1}{2j+\la} - \frac{1}{2k+\la}\Big].
$$
From the equation $\Gamma(t+1) = t\,\Gamma(t)$ it follows that
$\psi(t+1) = \psi(t)+1/t$, and the latter implies
\be \lb{m'}
    [\log m(\lambda)]'
 = - 2\sum_{i=2j}^{2k-1}\frac{1}{i+2\lambda}
   - \Big[\frac{1}{2j+\lambda} - \frac{1}{2k+\lambda}\Big],
\ee
whence
$$
    [\log m(\la)]''
 =  4\,\sum_{i=2j}^{2k-1}\frac{1}{(i+2\lambda)^2}
    + \Big[\frac{1}{(2j+\lambda)^2} - \frac{1}{(2k+\lambda)^2}\Big] > 0\,,
$$
and that proves the left-hand inequalities in \rf[b<]-\rf[b>].


2) We approach in the same way to the proof of the right-hand
inequalities in \rf[b<] and \rf[b>], by taking the
difference of the logarithms of the middle and the right-hand terms,
\be \lb{h}
    h(\la)
 := \log m(\la) - \log r(\la)
  = \log m(\la) - (2\la-2)\log\frac{j+\la}{k+\la}\,.
\ee We need to show that $h(\la) \ge 0$ for $\la \in [0,1]$ and that
$h(\la) < 0$ otherwise. Since $h(0)= h(1)=0$ by \rf[lmr], it
suffices to show that $h'(\la) < 0$ for $\la > 1$ and that $h''(\la)
< 0$ for $\la \in (-\frac{1}{2}, 1]$.

\smallskip
2a) Let us show that $h'(\la)\leq 0$ for $\la\ge 1$.
From \rf[h] using \rf[m'],  we obtain
\be \lb{h'}
    h'(\la)
 = - 2\sum_{i=2j}^{2k-1}\frac{1}{i+2\lambda}
   - \Big[\frac{1}{2j+\lambda} - \frac{1}{2k+\lambda}\Big]
   - 2\log\,\frac{j+\lambda}{k+\lambda}
   - (2\lambda-2)\Big[\frac{1}{j+\lambda} -
   \frac{1}{k+\lambda}\Big].
\ee
For the sum, since the function $f(x) = (x+2\la)^{-1}$ is decreasing,
we have
$$
   - 2\sum_{i=2j}^{2k-1}\frac{1}{i+2\la}
 < - 2 \int_{2j}^{2k} \frac{1}{x+2\la}\,dx
 =   2\,\log\frac{j+\la}{k+\la}\,,
$$
hence
\be \lb{h'1}
    h'(\la)
 <  - \Big[\frac{1}{2j+\la} - \frac{1}{2k+\la}\Big]
    - (2\la-2)\Big[\frac{1}{j+\la} - \frac{1}{k+\la}\Big],
\ee
and for $\la > 1$ and $j < k$, the right-hand side is negative.
Thus, $h'(\la) < 0$ for $\la > 1$.

\medskip
2b) Next, we prove that if $\la \in (-\frac{1}{2},1]$, then
$h''(\la)<0$. From \rf[h'], we derive
\ba
    h''(\la)
&=& 4\,\sum_{i=2j}^{2k-1}\frac{1}{(i+2\la)^2}
    + \Big[\frac{1}{(2j+\la)^2}-\frac{1}{(2k+\la)^2}\Big]
    \lb{h''}\\
&& \quad   - 4\, \Big[\frac{1}{j+\la} - \frac{1}{k+\la}\Big]
    + (2\la-2)\Big[\frac{1}{(j+\la)^2}-\frac{1}{(k+\la)^2}\Big]\,.
    \lb{h''2}
\ea
The first term in the right-hand side is estimated as follows
\baa
    4\,\sum_{i=2j}^{2k-1}\frac{1}{(i+2\la)^2}
&=&  4\,\sum_{i=2j+1}^{2k}\frac{1}{(i+2\la)^2}
    + \Big[\frac{1}{(j+\la)^2} - \frac{1}{(k+\la)^2}\Big] \\
&<&  2\, \Big[\frac{1}{j+\la} - \frac{1}{k+\la}\Big]
    + \Big[\frac{1}{(j+\la)^2} - \frac{1}{(k+\la)^2}\Big], \\
\eaa where for the sum we have used the inequality
$\sum_{i=2j+1}^{2k} (i+2\la)^{-2} < \int_{2j}^{2k}
(x+2\la)^{-2}\,dx$.

Next, for $\la \in (-\frac{1}{2},1]$ and $x\ge \Frac{1}{2}$ the function
$f(x)=(2x+\la)^{-2}-(x+\la)^{-2}$ is increasing, hence
for the second term in \rf[h''] we have
$$
   \Big[\frac{1}{(2j+\la)^2}-\frac{1}{(2k+\la)^2}\Big]
 < \Big[\frac{1}{(j+\la)^2}-\frac{1}{(k+\la)^2}\Big]\,.
$$
Substituting the above upper bounds in the expression \rf[h'']-\rf[h''2]
for $h''(\la)$, we obtain
\be \lb{h''3}
    h''(\la)
 <  - 2\, \Big[\frac{1}{j+\la} - \frac{1}{k+\la}\Big]
    + 2\la\,\Big[\frac{1}{(j+\la)^2}-\frac{1}{(k+\la)^2}\Big]
 =  - \frac{2(k-j)(kj-\la^2)}{(j+\la)^2(k+\la)^2}
 <  0\,,
\ee
since $1 \le j < k$ and $\la \in (-\frac{1}{2},1]$.
\qed


\begin{proposition}\lb{wbj/bk}
Let $j,k\in \N$, $j < k$. Then the coefficients $\wt\beta_k$ in
\rf[wb] satisfy the following relations.

\medskip
(i)\quad If $\,-\frac{1}{2} <\la\leq 0\,$ or $\,\lambda\ge 1$, then
\be \lb{wb<}
    \Big(\frac{j-\frac{1}{2}}{k-\frac{1}{2}}\Big)^{2\la-2}
\le \frac{\wt\beta_{k}^2}{\wt\beta_j^2}
\le \Big(\frac{j-\frac{1}{2}+\la}{k-\frac{1}{2}+\la}\Big)^{2\la-2}\,.
\ee

(ii)\quad If  $\,0<\lambda\le 1$,  then
\be \lb{wb>}
    \Big(\frac{j-\frac{1}{2}}{k-\frac{1}{2}}\Big)^{2\la-2}
\ge \frac{\wt\beta_{k}^2}{\wt\beta_j^2}
\ge \Big(\frac{j-\frac{1}{2}+\la}{k-\frac{1}{2}+\la}\Big)^{2\la-2}\,.
\ee
\end{proposition}

\proof
By equality \rf[wbb], we have
$$
     \wt\beta_j = \beta_{j-\frac{1}{2}}\,, \qquad
     \wt\beta_k = \beta_{k-\frac{1}{2}}\,.
$$
Then all the relations throughout \rf[b1]-\rf[h''3] remain valid
with the substitution
$$
     j \to j-\Frac{1}{2}, \qquad k \to k - \Frac{1}{2}\,.
$$
The only exception is inequality \rf[h''3] which fails for $j=1$,
$k=2$, and $\la \in [\frac{\sqrt{3}}{2},1]$, since the factor
$\big[(k-\frac{1}{2})(j-\frac{1}{2}) - \la^2\big]$ is not positive
then.

Let us prove that $\wt h(\la) \ge 0$ in this case as well. Since
$\wt h(1) = 0$, it is sufficient to prove that $\wt h'(\la) < 0$ for
$\la \in [\frac{\sqrt{3}}{2},1]$ and $j = 1$, $k=2$. We have
$$
   \wt h'(\la)\Big|_{j,k} = h'(\la)\Big|_{j-\frac{1}{2},k-\frac{1}{2}}
$$
so substituting $j = \frac{1}{2}$, $k = \frac{3}{2}$ into \rf[h'1],
we find that for $\la \in [\frac{3}{4},1] \supset [\frac{\sqrt{3}}{2},1]$
\baa
    \wt h'(\la)\Big|_{1,2}
 = h'(\la)\Big|_{\frac{1}{2},\frac{3}{2}}
&<& - \Big[\frac{1}{1+\la} - \frac{1}{3+\la}\Big]
    - (2\la-2)\Big[\frac{1}{\frac{1}{2} + \la}
          - \frac{1}{\frac{3}{2} + \la}\Big] \\
&\le& - \Big[\frac{1}{1+\la} - \frac{1}{3+\la}\Big]
      + \Frac{1}{2} \Big[\frac{1}{\frac{1}{2} + \la}
          - \frac{1}{\frac{3}{2} + \la}\Big] \\
& = &   -\, \frac{2}{(1+\la)(3+\la)} + \frac{2}{(1+2\la)(3+2\la)} < 0.
\eaa


\section{Three lemmas} \lb{lemmas}

In the next two sections, we deal with lower and upper estimates
for the sums $\sum_{j=1}^\l f(j)$, in particular for $f=F_\nu$, where
$F_1, F_2$ are given in \rf[F] below.
For that purpose, we need the following three lemmas.

We use the following notation:
$$
   \sum_{i=1}^{\l}\!{}^{''} f(i)
= \Frac{1}{2}f(1) + \sum_{i=2}^{\l-1} f(i) + \Frac{1}{2} f(\l)\,.
$$

\begin{lemma} \lb{leint0}
For a convex integrand $f$, we have
$$
    \sum_{i=1}^{\l} f(i)
\le \int_{\frac{1}{2}}^{\l+\frac{1}{2}} f(x)\,dx,\qquad
    \sum_{i=1}^{\l}\!{}^{''} f(i)
\ge \int_1^{\l} f(x)\,dx\,.
$$
\end{lemma}

\proof The inequalities reveal well-known properties of the midpoint
and the trapezoidal quadrature formulas relative to the
corresponding integrals. \qed

\begin{lemma} \lb{f}
For $\la > - \frac{1}{2}$, the functions
\be \lb{F}
   F_1(x) = x^{2\la} (x + \Frac{\la}{2})^2 (x + \la)^2, \qquad
   F_2(x) = x^2 (x + \Frac{\la}{2})^2 (x + \la)^{2\la}
\ee
are convex on $[\frac{1}{2},\infty)$ and increasing on $[1,\infty)$.
\end{lemma}

\proof
1) For $\la \ge 1$, all the factors of $F_1$, $F_2$ in \rf[F]
are convex, positive and increasing on $[0,\infty)$, hence the statement.

2) For $\la \in [0,1]$ the functions
$$
   u_1(x) := x^\la (x+\la), \qquad u_2(x) := x (x+\la)^\la
$$
are non-negative and increasing on $[0,\infty)$. Further, $u_2$ is
convex on $[0,\infty)$, because it can be written in the form
$$
     u_2(x) = (x+\la)^{\la+1} - \la (x+\la)^\la\,,
$$
where both terms are convex for $\la \in [0,1]$, whereas $u_1$ is
convex on $[\frac{1}{2},\infty]$ because
$$
   u_1''(x)
 = [x^{\la+1} + \la x^\la]''
 = \la x^{\la-2}\Big[(\la+1)x + \la(\la-1)\Big]
 > \la x^{\la-2}\Big[x - \Frac{1}{4}\Big]
\ge 0, \quad x \ge \Frac{1}{2}\,.
$$
Therefore, both $F_1(x) = [u_1(x)]^2 (x+\frac{\la}{2})^2$ and
$F_2(x) = [u_2(x)]^2 (x+\frac{\la}{2})^2$ are convex on
$[\frac{1}{2},\infty)$ and increasing on $(0,\infty)$.

3) Let $\la \in (-\frac{1}{2},0]$. Then
$$
   u_1'(x)
 = x^{\la-1}\Big[(\la+1) x + \la^2\Big]
 > 0, \quad x > 0,
$$
and
$$
   u_2'(x)
 = (x+\la)^{\la-1}\Big[(\la+1) (x+\la) - \la^2\Big]
\ge (x+\la)^{\la-1}\Big[(\la+1)^2 - \la^2\Big]
 > 0, \quad x \ge 1,
$$
hence $F_1$ and $F_2$ are increasing on $[1,\infty)$.
Further, the function
$$
   v_1(x)
:= x^\la(x+\Frac{\la}{2})(x+\la)
 = x^{\la+2} + \Frac{3\la}{2}\,x^{\la-1} + \Frac{\la^2}{2} x^\la
$$
is convex for $x > 0$ because all the terms are convex
for $\la \in (-\frac{1}{2},0]$,
hence $F_1(x) = [v_1(x)]^2$ is convex whenever $v_1$ is nonnegative, i.e.,
for $x > -\la$, thus for $x \ge \frac{1}{2}$.
Finally, for
$$
   v_2(x)
:= x(x+\Frac{\la}{2})(x+\la)^\la
 = y^{\la+2} - \Frac{3\la}{2}\,y^{\la+1}
    + \Frac{\la^2}{2} y^\la\,,\quad y = x+\la,
$$
we obtain
$$
   v_2''(x)
 = y^{\la-2}\Big[(\la+2)(\la+1)y^2 - \Frac{3}{2}\la^2(\la+1)y
    + \Frac{1}{2}\la^3(\la-1)\Big] =: y^{\la-2} p_2(y)\,,
$$
and it is easy to check that, for $\la \in (-\frac{1}{2},0]$,
the quadratic polynomial $p_2$ has no real zeros. Hence, $v_2$
is convex and so is $F_2(x) = [v_2(x)]^2$ for $x \ge \frac{1}{2}$.
\qed


\begin{lemma} \lb{leint}
Let $\,\al_i>0$, $\gamma_{\min} \le \gamma_i \le \gamma_{\max}$,
$1\le i \le r$, and let
$$
    f(x) := (x+\gamma_1)^{\al_1}(x+\gamma_2)^{\al_2}
              \cdots(x+\gamma_r)^{\al_r}, \qquad
    s := \sum_{i=1}^{r}\al_i\,.
$$
Then, for any $\,x > x_0$,  where $x_0 + \gamma_{\min} \ge 0$,
we have
\be \lb{int}
   \frac{1}{s+1}\, \Big[(t+\gamma_{\min})f(t)\Big]_{x_0}^x
 < \int\limits_{x_0}^{x} f(t)\,dt
 < \frac{1}{s+1}\,(x+\gamma_{\max})f(x)\,.
\ee
\end{lemma}

\proof
Set
$$
   G(x) := \frac{1}{s+1}\,(x+\gamma_{\min}) f(x), \qquad
   F(x) := \frac{1}{s+1}\,(x+\gamma_{\max}) f(x)\,.
$$
It suffices to show that $G'(t) < f(t) < F'(t)$ for $x_0 \le t \le x$.
We have
$$
   G'(t)
=  \frac{1}{s+1}
   \Big[1 + \sum_{i=1}^r \alpha_i \frac{t+\gamma_{\min}}{t+\gamma_i}\Big]f(t)
\le \frac{1}{s+1}
   \Big[1 + \sum_{i=1}^r \alpha_i \Big] f(t) = f(t)\,,
$$
and similarly
$$
   F'(t)
=  \frac{1}{s+1}f(t)
   \Big[1 + \sum_{i=1}^s \alpha_i \frac{t+\gamma_{\max}}{t+\gamma_i}\Big]
\ge \frac{1}{s+1}
   \Big[1 + \sum_{i=1}^s \alpha_i \Big] f(t) = f(t)\,.\qquad \Box
$$

\begin{remark} \rm
We can refine the upper estimate as follows:
$$
   \int_{x_0}^{x} f(t)\,dt
 < \frac{1}{s+1} [f(x)]^\frac{s+1}{s}\,.
$$
Indeed, with $F(x) := \frac{1}{s+1} [f(x)]^\frac{s+1}{s}$, it
suffices to show that $\,F'(t)\ge f(t)\,$ for every $\,t > x_0$. We
have the equivalent relations
$$
    F'(t) = \frac{1}{s}\,\big[f(t)\big]^{\frac{1}{s}}\,f'(t)\,
\ge f(t)
\lr \big[f(t)\big]^{\frac{1}{s}}
\ge \frac{s}{\frac{f'(t)}{f(t)}}\,,
$$
and the latter is simply the inequality between the
geometric and harmonic means
$$
    \Big(\prod(x+\gamma_i)^{\alpha_i}\Big)^{\frac{1}{\sum\alpha_i}}
\ge \frac{\sum\alpha_i}{\sum\frac{\alpha_i}{x+\gamma_i}}\,.
$$
\end{remark}

\section{An upper bound for $\|\A_m\|_\infty$ for $\la > 2$}
\lb{infty}

\begin{proposition}\lb{A0}
For $\la>2$, we have \be \lb{A1}
      \|\A_m\|_\infty
\le \frac{4}{(\la + 2)(\la + 3)}\,
      m(m + \Frac{\la}{2})(m + \la)(m + \Frac{3\la}{2} +3)\,.
\ee
\end{proposition}

\proof Let us recall that
$$
   \|\A_m\|_\infty = \max_k \sum_j |a_{kj}|\,,
$$
and, as is seen from  \rf[A_m], $\,a_{k,j}>0$.

For a fixed $\,k$, $\,1\le k\le m$, we consider the sum of the
elements in the $k$-th row of $\A_m$,
$$
    \sum_{j=1}^{m} a_{kj}
 =  \sum_{j=1}^{k-1}\frac{\beta_k}{\beta_j}\,a_{jj}
   + a_{kk}
+ \sum_{j=k+1}^{m}\frac{\beta_j}{\beta_{k}}\,a_{kk}.
$$
By \rf[al] and by \rf[b<],
\be \lb{ab}
   a_{jj}
=  c_0 f_0(j), \qquad
   \frac{\beta_k}{\beta_j} \le \Big(\frac{j+\la}{k+\la}\Big)^{\la-1},
   \quad j < k, \quad \la > 0\,,
\ee
where
$$
c_0 := \frac{4}{2\la+1}\,,\qquad
    f_0(x) := x(x+\Frac{\la}{2})(x+\la)\,,
$$
hence
\be \lb{a0}
    \sum_{j=1}^{m} a_{kj}
\le c_0
    \Big[\sum_{j=1}^{k-1}f_0(j)\Big(\frac{j+\la}{k+\la}\Big)^{\la-1}
    + f_0(k)
    + f_0(k) \sum_{j=k+1}^{m}\Big(\frac{k+\la}{j+\la}\Big)^{\la-1}
    \Big]\,.
\ee For the first sum, since $f(x) = f_0(x)(x+\la)^{\la-1}$ is
increasing, we apply an integral estimate and then Lemma \ref{leint}
to obtain
$$
      \sum_{j=1}^{k-1}f_0(j)\Big(\frac{j+\la}{k+\la}\Big)^{\la-1}
\le \int_{1}^{k} f_0(x)\Big(\frac{x+\la}{k+\la}\Big)^{\la-1}\,dx
\le  \frac{1}{\lambda+3}\,(k+\la)f_0(k) \,.
$$
For the second sum, since $g(x) = 1/(x+\la)^{\la-1}$ is decreasing
(and $\la>2$), an integral estimate gives
$$
    f_0(k) \sum_{j=k+1}^{m} \Big(\frac{k+\la}{j+\la}\Big)^{\la-1}
\le f_0(k) \int\limits_{k}^{m} \Big(\frac{k+\la}{x+\la}\Big)^{\la-1}\,dx
 =  \frac{1}{\la-2} (k+\la)f_0(k)
    \Big[1 - \Big(\frac{k+\la}{m+\la}\Big)^{\la-2} \Big]\,.
$$
Replacement in the right-hand of \rf[a0] yields
\be \lb{a1}
     \sum_{j=1}^{m}a_{kj}
  <  c_0\, (k+\la) f_0(k)
     \Bigg[\frac{1}{\la+3}+\frac{1}{\la-2}-\frac{1}{\la-2}
     \Big(\frac{k+\la}{m+\la}\Big)^{\la-2}\Bigg]
     + c_0 f_0(k)
     =: A + B\,.
\ee

1) We estimate $A$ as follows.
\ba
     A
 &=& \frac{4}{2\la+1}\, (k+\la)f_0(k)
     \Bigg[\frac{1}{\la+3}+\frac{1}{\la-2}-\frac{1}{\la-2}
     \Big(\frac{k+\la}{m+\la}\Big)^{\la-2}\Bigg] \nonumber \\
 &=& \frac{4}{2\la+1}\, (k+\la)f_0(k)
     \Bigg[\frac{2\la+1}{(\la+3)(\la-2)} - \frac{1}{\la-2}
     \Big(\frac{k+\la}{m+\la}\Big)^{\la-2}\Bigg] \nonumber \\
&=& \frac{4(m+\la)^4}{(\la+3)(\la-2)}\,
     \frac{f_0(k)}{(k+\la)^3}\;
     \Big(\frac{k+\la}{m+\la}\Big)^4
     \Bigg[1-\frac{\la+3}{2\la+1}
     \Big(\frac{k+\la}{m+\la}\Big)^{\la-2}\Bigg] \nonumber \\
 &=& \frac{4(m+\la)^4}{(\la+3)(\la-2)}\,
     \psi_\la(k)\, \phi_{\la}(y)\,,
     \lb{A}
\ea
where in the last line we set
$$
    \phi_{\la}(y)
 := y^4  - \frac{\la+3}{2\la+1}\,y^{\la+2}\,,\qquad
    y = \frac{k+\la}{m+\la}\in [0,1]\,, \qquad
    \psi_{\la}(k) := \frac{f_0(k)}{(k+\la)^3}\,.
$$
Let us evaluate $\phi_\la(y)$ and $\psi_\la(k)$.
On $[0,1]$, for a fixed $\la > 2$, the function $\,\phi_\la\,$ has a unique
local extremum, a maximum, which is attained at
$$
    y_\la
 =  \Big(\frac{4(2\la+1)}{(\la+2)(\la+3)}
      \Big)^{\frac{1}{\la-2}}
 = \Big(1 - \frac{(\la-1)(\la-2)} {(\la+2)(\la+3)}
   \Big)^{\frac{1}{\la-2}}
\in (0,1)\,.
$$
Then
\be \lb{phi}
    \phi_{\la}(y)
\le \phi_{\la}(y_\la)
 =  \frac{\la - 2}{\la + 2}\,y_\la^4
 < \frac{\la - 2}{\la + 2}\,.
\ee
The function $\psi_\la(x) = \frac{x(x+\frac{\la}{2})(x+\la)}{(x+\la)^3}$
is increasing (since $h(x) = \frac{x+a}{x+b}$ is increasing for $a < b$),
thus
\be \lb{psi}
   \psi_\la(k) = \frac{f_0(k)}{(k+\la)^3} \le \frac{f_0(m)}{(m+\la)^3}\,.
\ee
Consequently, putting the estimates \rf[phi]-\rf[psi] into \rf[A], we obtain
$$
     A \le \frac{4}{(\la+2)(\la+3)}\, (m+\la)f_0(m)\,.
$$

2) For $B$ in \rf[a1] we use the trivial upper estimate
$$
     B = c_0 f_0(k) \le c_0 f_0(m)
 = \frac{4}{2\la+1}\,f_0(m)\,.
$$

3) Thus, from \rf[a1], we derive
\baa
      \sum_{k=j}^{m}a_{kj} \le A + B
&\le& \frac{4}{(\la+2)(\la+3)}\,
      \Big(m+\la + \frac{(\la+2)(\la+3)}{2\la+1} \Big)f_0(m) \\
&\le& \frac{4}{(\la+2)(\la+3)}\,
      (m+\Frac{3\la}{2}+3) f_0(m)\,,
\eaa where we have used that $\frac{(\la+2)(\la+3)}{2\la+1} <
\Frac{\la}{2} + 3$ for $\la > 2$. Hence,
\begin{equation} \lb{A1'}
     \|\A_m\|_{\infty}
 \le \frac{4}{(\la+2)(\la+3)}\,m(m+\Frac{\la}{2})
     (m+\la)(m+\Frac{3\la}{2}+3) \,,
\end{equation}
and \rf[A1] is proved. \qed

\begin{proposition} \lb{pB1}
For $\la>2$, and $n\in\N$, we have
\ba
      \|\B_n\|_\infty
&\le& \frac{1}{(\la + 2)(\la + 3)}\, n(n+\la)(n+2\la)(n+3\la+6)
         \lb{B1'} \\
&\le& \frac{1}{(\la + 2)(\la + 3)}\, n(n+2\la+ 2)^3\,. \lb{B1}
\ea
\end{proposition}

\proof
Recall that
\be \lb{AB}
   \B_n = \left\{\begin{array}{ll}
   4 \A_m, & n = 2m  \\
   4 \wt\A_m, & n = 2m-1\,.
   \end{array} \right.
\ee Let us rewrite \rf[A1'] as \be \lb{AK}
   \|\A_m\|_\infty
\le 4 c_1 K_1(m) = 4 c_1 K_1(\Frac{n}{2}), \qquad n = 2m\,,
\ee
where
$$
   c_1 := \frac{1}{(\la+2)(\la+3)}\,, \qquad
   K_1(m) :=  m(m + \Frac{\la}{2})(m + \la)(m + \Frac{3\la}{2} +3)\,.
$$
We derived this upper bound from two estimates in \rf[ab], namely
$$
   a_{jj}
=  c_0\,f_0(k), \qquad
   \frac{\beta_k}{\beta_j} \le \Big(\frac{j+\la}{k+\la}\Big)^{\la-1}, \quad
   j < k\,.
$$
Now, we note that, by \rf[wb<] and \rf[wtal],  we have similar estimates
$$
   \wt a_{jj}
=  c_0\,f_0(k-\Frac{1}{2}), \qquad
   \frac{\wt\beta_k}{\wt\beta_j}
   \le \Big(\frac{j-\frac{1}{2}+\la}{k-\frac{1}{2}+\la}\Big)^{\la-1}, \quad
   j < k\,,
$$
and it is easy to see that all the inequalities for the sum
$\sum_j \wt a_{kj}$ throughout \rf[a0]-\rf[A1']
remain valid with
the substitution $(j,k)\to (j-\frac{1}{2},k-\frac{1}{2})$,
hence
\be \lb{wtAK}
   \|\wt\A_m\|_\infty
\le 4 c_1 K_1(m-\Frac{1}{2}) = 4 c_1 K_1(\Frac{n}{2}),\qquad n = 2m-1.
\ee
Now, from \rf[AB], \rf[AK] and \rf[wtAK], we obtain that for any $n\in \N$
\baa
      \|\B_n\|_\infty
\le 16 c_1 K_1\big(\Frac{n}{2}\big)
&=& c_1 n(n+\la)(n+2\la)(n+3\la+6) \\
&\le& c_1 n(n+2\la+2)^3,
\eaa
where the last inequality follows by relation
between geometric and arithmetic means, namely $abc \le (\frac{a+b+c}{3})^3$,
with $(a,b,c) = (n +\la, n+2\la, n+3\la+6)$. This proves \rf[B1']-\rf[B1].
\qed


\section{Lower and upper estimates for $\|\A_m\|_F$
for $\la > - \frac{1}{2}$} \lb{Fro}

\begin{proposition}
For $\la > - \frac{1}{2}$, we have \be \lb{AF}
        c_4\, (m+\la')(m+\la'+4) F_0(m)
\;\le\; \|\A_m\|_F^2
\;\le\; c_4\, (m+\la+\la'' +\Frac{5}{2})^2 F_0(m+\Frac{1}{2}),
\ee
where $\la' := \min\,\{0,\la\}$, $\la'' := \max\,\{0,\la\}$, and
$$
    c_4 := \frac{4}{(2\la+1)^2(2\la+5)}\,, \qquad
    F_0(x) := [f_0(x)]^2 = x^2 (x+\Frac{\la}{2})^2 (x+\la)^2\,.
$$
\end{proposition}

\proof By the definition of the Frobenius norm,
$$
   \|\A_m\|_F^2 := \sum_{j,k=1}^m {a_{kj}^2}\,.
$$
Since matrices $\{\A_m\}$ are symmetric and embedded, we have \be
\lb{N0}
      \|\A_k\|_F^2 - \|\A_{k-1}\|_F^2
\;=\; 2 \sum_{j=1}^k\!{}^{'} a_{kj}^2
\;=\; 2 \sum_{j=1}^k\!{}^{'} \frac{\beta_k^2}{\beta_j^2}\,
       a_{jj}^2\,.
\ee
where $\sum'$ means that the last summand is halved.
Recall that by \rf[al]
$$
   a_{jj}^2
 = c_0^2\, [f_0(j)]^2 =: c_2 F_0(j),
   \qquad c_2 := \frac{16}{(2\la+1)^2}\,.
$$

1) {\it The case $\la \in (-\frac{1}{2},0] \cup [1,\infty)$}.
In that case, by \rf[b<],
$$
      \Big(\frac{j}{k}\Big)^{2\la-2}
\le   \frac{\beta_k^2}{\beta_j^2}
\le \Big(\frac{j+\la}{k+\la}\Big)^{2\la-2},
$$
so we obtain from \rf[N0]
\be \lb{N}
     2\,c_2 \sum_{j=1}^k\!{}^{'} f_1(j)
\le  \|\A_k\|_F^2 - \|\A_{k-1}\|_F^2
\le  2\,c_2 \sum_{j=1}^k\!{}^{'} f_2(j)\,,
\ee
where
\ba
    f_1(x)
&:=& F_0(x) \Big(\frac{x}{k}\Big)^{2\la-2}
 =  \frac{1}{k^{2\la-2}}\,x^{2\la}
            (x+\Frac{\la}{2})^2(x+\la)^2, \lb{f1} \\
    f_2(x)
&:=& F_0(x) \Big(\frac{x+\la}{k+\la}\Big)^{2\la-2}
 =  \frac{1}{(k+\la)^{2\la-2}}\,
            x^2(x+\Frac{\la}{2})^2(x+\la)^{2\la}\,. \lb{f2}
\ea Note that, by Lemma \ref{f}, both functions are convex on
$\,[\frac{1}{2},\infty)\,$ and monotonely increasing on
$[1,\infty)$, and that
$$
   f_1(k) = f_2(k) = F_0(k)\,.
$$
Set
$$
    \la' := \min\, \{0,\Frac{\la}{2},\la\} = \min\,\{0,\la\}\,, \qquad
    \la'':= \max\, \{0,\Frac{\la}{2},\la\} = \max\,\{0,\la\}\,.
$$
Those will play the roles of $\gamma_{\max}$ and $\gamma_{\min}$
when we apply Lemma \ref{leint}.

1a) For the upper estimate, since $f_2$ is convex and increasing, we
have by Lemmas \ref{leint0} and \ref{leint} for $k\geq 2$, \baa
    \sum_{j=1}^{k-1} f_2(j)
\le \int_{\frac{1}{2}}^{k-\frac{1}{2}} f_2(x) dx \le
\frac{k-\Frac{1}{2}+\la''}{2\la+5} f_2(k-\Frac{1}{2}) \le
\frac{k-\Frac{1}{2}+\la''}{2\la+5}  f_2(k)\,, \eaa so that, for
$k\geq 1$,
$$
    \sum_{j=1}^k\!{}^{'} f_2(j)
\le \frac{k-\Frac{1}{2}+\la''}{2\la+5} f_2(k) + \frac{1}{2} f_2(k)
    = c_3\,(k+\la''+\la+2) f_2(k), \qquad c_3 := \frac{1}{2\la+5}\,,
$$
hence
\be \lb{Ng2}
    \|\A_k\|_F^2 - \|\A_{k-1}\|_F^2
\le 2\,c_2 c_3\,(k+\la''+\la+2) f_2(k)
 =  2\,c_2 c_3\,(k+\la''+\la+2) F_0(k)\,.
\ee
Then,
$$
    \|\A_m\|_F^2
  = \sum_{k=1}^m  \Big(\|\A_k\|_F^2 - \|\A_{k-1}\|_F^2\Big)
\le 2\,c_2 c_3 \sum_{k=1}^m g_2(k)\,,
$$
where $g_2(x) = (x+\la''+\la +2) F_0(x)$ is convex, and by Lemmas
\ref{leint0} and \ref{leint} we obtain
\be \lb{up}
      \sum_{k=1}^m g_2(k)
\le \int_{\frac{1}{2}}^{m+\frac{1}{2}} g_2(x)\,dx
 \le  \frac{1}{8}\, (m+\la+\la''+\Frac{5}{2})\,g_2(m+\Frac{1}{2})
  =  \frac{1}{8} (m+\la+\la'' +\Frac{5}{2})^2 F_0(m+\Frac{1}{2}) \,,
\ee and this proves the upper estimates in \rf[AF] for  $\,\la \in
(-\frac{1}{2},0] \cup [1,\infty)$, with the constant $c_4 =
\frac{1}{4}c_2 c_3$.

\medskip
1b) For the lower estimate, we get by Lemmas \ref{leint0} and \ref{leint},
\baa
   \sum_{j=1}^k\!{}^{'} f_1(j)
&=&  \frac{1}{2}f_1(1) + \sum_{j=1}^k\!{}^{''} f_1(j)
\;\ge\;\frac{1}{2}f_1(1) + \int_1^k f_1(x)\,dx \\
&\ge& \frac{1}{2}f_1(1) + \frac{k+\la'}{2\la+5}  f_1(k)
    - \frac{1+\la'}{2\la+5} f_1(1)
\;\ge\;  \frac{k+\la'}{2\la+5} f_1(k)\,,
\eaa
hence
\be \lb{Ng1}
    \|\A_k\|_F^2 - \|\A_{k-1}\|_F^2
\ge 2 c_2 c_3 (k+\la') f_1(k)
 =  2 c_2 c_3 (k+\la') F_0(k)\,.
\ee
Then,
$$
    \|\A_m\|_F^2 \ge 2 c_2 c_3 \sum_{k=1}^m g_1(k)\,,
$$
where $g_1(x) = (x+\la') F_0(x)$ is convex, therefore, by Lemmas
\ref{leint0} and \ref{leint},
\ba
      \sum_{k=1}^m g_1(k)
&=&  \frac{1}{2}g_1(1) + \sum_{j=1}^k\!{}^{''} g_1(j) + \frac{1}{2}
g_1(m)
 \ge \frac{1}{2} g_1(1) + \int_1^m g_1(x)\,dx + \frac{1}{2} g_1(m)
     \nonumber \\
&\ge&  \frac{1}{2} g_1(1) + \frac{m+\la'}{8} g_1(m)
      - \frac{1+\la'}{8} g_1(1) + \frac{1}{2} g_1(m)
\ge \frac{m+\la'}{8} g_1(m) + \frac{1}{2} g_1(m) \nonumber \\
&\ge& \frac{1}{8} (m + \la'+ 4)\, g_1(m)
 = \frac{1}{8} (m+\la')(m +\la' + 4)\, F_0(m)\,, \lb{low}
\ea
and the lower estimate in \rf[AF] follows, with $c_4 = \frac{1}{4}c_2c_3$.

\medskip
2) {\it The case $\la \in [0,1]$.} In that case, by \rf[wb<], we have
$$
      \Big(\frac{j+\la}{k+\la}\Big)^{2\la-2}
\le   \frac{\beta_k^2}{\beta_j^2}
\le \Big(\frac{j}{k}\Big)^{2\la-2},
$$
so we obtain
$$
     2 c_2 \sum_{j=1}^k\!{}^{'} f_2(j)
\le  \|\A_k\|_F^2 - \|\A_{k-1}\|_F^2  
\le  2 c_2 \sum_{j=1}^k\!{}^{'} f_1(j)\,,
$$
i.e., the same inequality as in \rf[N], but with $f_1$ and $f_2$
interchanged.

\medskip
2a) Then the upper estimates will run in the same way
only with $f_1$ instead of $f_2$, and because
\be \lb{ff1}
   f_1(k) = f_2(k) = F_0(k)
\ee
we arrive at the same inequality \rf[Ng2], so that the final upper
estimate for $\|\A_m\|_F^2$ for $\la \in [0,1]$ is the same as \rf[up].

2b) Similarly, the lower estimates for $\la \in [0,1]$
will run in the same way only with $f_2$ instead of  $f_1$,
and because of \rf[ff1] we arrive at the same inequality
\rf[Ng1], so that the final lower estimate for $\|\A_m\|_F^2$
for $\la \in [0,1]$ is also the same as \rf[low].
\qed


\begin{proposition} \lb{pB2}
For $n \in \N$ and $\la > - \frac{1}{2}$, we have
\ba
    c_5^2\, (n+8)n^3(n+\la)^2(n+2\la)^2 \;
\le &\|\B_n\|_F^2&
\le \; c_5^2\, (n+2\la+2)^8,
    \quad \la \ge 0; \lb{B2} \\[0.5ex]
    c_5^2\, (n+2\la+8)n^2(n+\la)^2(n+2\la)^3 \;
\le &\|\B_n\|_F^2&
\le \; c_5^2\, (n+\la+2)^8,
    \quad\; \la \in (-\Frac{1}{2},0] \lb{B3}
\ea
where
$$
    c_5^2 := \frac{1}{16}\,c_4 = \frac{1}{4(2\la+1)^2(2\la+5)}\,.
$$
\end{proposition}

\proof
Recall again that
\be \lb{AB1}
   \B_n = \left\{\begin{array}{ll}
   4\A_m, & n = 2m; \\
   4\wt\A_m, & n = 2m-1,
   \end{array} \right.
\ee
and rewrite \rf[AF] as
\be \lb{AFK}
   c_4 K_{2,\la}(m) \le \|\A_m\|_F^2 \le c_4 K_{3,\la}(m),
   \qquad n = 2m\,.
\ee
Then, for odd $n =2m-1$, by the same arguments as in the proof of Proposition
\ref{pB1}, we obtain
\be \lb{AFK1}
    c_4 K_{2,\la}(m-\Frac{1}{2})
\le \|\wt\A_m\|_F^2
\le c_4 K_{3,\la}(m-\Frac{1}{2}), \qquad n =2m-1\,,
\ee
so that, for all $n\in \N$,
\be
    16 c_4 K_{2,\la}\Big(\frac{n}{2}\Big)
\le \|\B_n\|_F^2
\le 16 c_4 K_{3,\la}\Big(\frac{n}{2}\Big)\,.
\ee

Simplifying $K_{2,\la}(\frac{n}{2})$ we obtain
\begin{eqnarray*}
    K_{2,\la}\Big(\frac{n}{2}\Big)
&=& \frac{1}{2^8} n^2 (n+\la)^2 (n+2\la)^2 (n+2\la')(n+2\la'+8)\,,
\end{eqnarray*}
and this gives the lower bounds in \rf[B2]-\rf[B3] with the constant
$$
    c_5^2 = \frac{16}{2^8}\,c_4 = \frac{1}{16}\,c_4\,.
$$

For the upper bounds we get \baa
    K_{3,\la}\Big(\frac{n}{2}\Big)
&=& \frac{1}{2^8} (n+1)^2 (n+\la+1)^2 (n+2\la+1)^2 (n+2\la+2\la''+5)^2 \\
&\le& \frac{1}{2^8} (n+\Frac{5}{4}\la+\Frac{1}{2}\la''+2)^8\,, \eaa
where we used the inequality $abcd \le (\frac{a+b+c+d}{4})^4$. The
last term does not exceed $\,2^{-8}(n+2\lambda+2)^8$, if $\,\la\geq
0$, and $\,2^{-8}(n+\lambda+2)^8$, if $\la\in (-\frac{1}{2},0)$.

That proves the upper bounds in \rf[B2]-\rf[B3]. \qed


\section{Proof of the main results} \lb{main}

Firstly, we will prove Theorem \ref{thm1} by establishing separately
the lower and the upper bounds therein.

\begin{theorem}
For the upper bounds, we have
\be \lb{c<}
   [c_n(\la)]^2
\le \left\{ \begin{array}{ll}
   \frac{1}{(\la+2)(\la+3)}\, n(n+2\la+2)^3, & \la > 2; \\
   \frac{1}{2(2\la+1)\sqrt{2\la+5}}\,(n+\la+\la''+2)^4,
       & \la > -\Frac{1}{2}. \\
   \end{array} \right.
\ee
where $\la'' = \max\,\{0,\la\}$.
\end{theorem}

\proof
We proved in Propositions \ref{pB1} and \ref{pB2} that
$$
  \|\B_n\|_\infty \le L_1(n,\la), \quad \la > 2, \qquad
  \|\B_n\|_F \le L_2(n,\la), \quad\la > -\frac{1}{2},
$$
where $L_\nu$ is the $\nu$-th line in \rf[c<],
and since $[c_n(\la)]^2 = \mu_{\max}(\B_n)$, and the
largest eigenvalue $\mu_{\max}(\B_n)$ is smaller than any matrix norm,
the upper bounds \rf[c<] follow.
\qed

\begin{theorem}
For the lower bounds, we have
\be \lb{c>}
    [c_n(\la)]^2
\ge \left\{ \begin{array}{ll}
   \frac{1}{4(\la+1)(\la+2)}\,n^2(n+\la)^2, & \la > 2; \\
   \frac{1}{(2\la+1)(2\la+5)}\,(n+\la)^2(n+2\la')^2,
      & \la > -\frac{1}{2}\,,
   \end{array} \right.
\ee
where $\la' = \min\,\{0,\la\}$.
\end{theorem}

\proof
1) The first inequality in \rf[c>] follows from second, since
$$
      \frac{1}{4(\la+1)(\la+2)} < \frac{1}{(2\la+1)(2\la+5)},
      \qquad n+2\la'= n \quad(\la > 0)\,.
$$

2) Let us prove the second inequality in \rf[c>]
splitting the cases $\la > 0$ and $-\frac{1}{2} < \la \le 0$.
We proved in Proposition \ref{pB2} that
\be \lb{B>}
   \|\B_n\|_F^2
\ge \left\{ \begin{array}{ll}
   c_5^2 n^3(n+8)(n+\la)^2(n+2\la)^2,  & \la \ge 0; \\
   c_5^2 n^2(n+\la)^2 (n+2\la)^3(n+2\la+8), & \la\in (-\frac{1}{2},0]\,,
   \end{array} \right.
\ee
where
$$
   c_5^2 = \frac{1}{4(2\la+1)^2(2\la+5)}\,.
$$
Next, we will need an expression for the trace of $\B_n$, which we
obtained in \cite[p. 17]{ans16}, \be \lb{tr}
   \tr (\B_n)
 = \left\{ \begin{array}{ll}
   c_6\,n(n+2)(n+2\la)(n+2\la+2)\,, & n = 2m; \\
   c_6 \Big[[(n+1)(n+2\la+1)]^2 - 2[(n+1)(n+2\la+1)]\Big]\,, & n = 2m-1\,,
   \end{array} \right.
\ee
where
$$
   c_6 = \frac{1}{4(2\la+1)}\,.
$$
From \rf[tr] we can get a common upper bound for both odd and even
$n$ as follows. For odd $n$, we obtain from \rf[tr] \ba
   \tr (\B_n)
&< & c_6 \big[(n+1)^2(n+2\la+1)^2 - (n+1)^2\big] \nonumber \\
&= & c_6 (n+1)^2 (n+2\la)(n+2\la+2), \qquad \la \ge 0, \lb{tr1} \ea
and \ba
   \tr (\B_n)
&\le & c_6\, [(n+1)^2(n+2\la+1)^2] - (n+2\la+1)^2 \nonumber \\
&=& c_6\, (n+2\la+1)^2 n(n+2), \qquad \la \in (-\Frac{1}{2},0]\,,
\lb{tr2} \ea and it is clear the both estimates \rf[tr1]-\rf[tr2]
give upper bounds for $\tr(\B_n)$ for even $n=2m$ in \rf[tr] as
well.

Set
$$
    c_7 = \frac{c_5^2}{c_6} = \frac{1}{(2\la+1)(2\la+5)}\,.
$$

2a) Then, for $\la \ge 0$, from \rf[mu>], \rf[B>] and \rf[tr1] we
have \baa
      \mu_{\max}(\B_n)
\ge \frac{\|\B_n\|_F^2}{\tr(\B_n)}
&\ge& c_7 \frac{n^3(n+8)(n+\la)^2(n+2\la)^2}{(n+1)^2(n+2\la)(n+2\la+2)} \\
&=:& c_7 n^2(n+\la)^2 \phi_\la(n) \\
&>&  c_7 n^2(n+\la)^2\,, \eaa since for $\la \ge 0$ and $n \ge 3$
$$
    \phi_\la(n)
:=  \frac{n(n+8)}{(n+1)^2}\frac{n+2\la}{n+2\la+2}
\ge \frac{n(n+8)}{(n+1)^2}\frac{n}{n+2}
\ge 1\,.
$$

2b) Similarly, for $\la \in(-\frac{1}{2},0]$,  from \rf[mu>],
\rf[B>] and \rf[tr2], we have \baa
      \mu_{\max}(\B_n)
 \ge  \frac{\|\B_n\|_F^2}{\tr(\B_n)}
&\ge& c_7 \frac{n^2(n+\la)^2 (n+2\la)^3(n+2\la+8)}{n(n+2)(n+2\la+1)^2} \\
& = & c_7 (n+\la)^2(n+2\la)^2 \psi_\la(n) \\
& > & c_7 (n+\la)^2(n+2\la)^2\,, \eaa since for $\la \in
(-\frac{1}{2}, 0]$ and $n \ge 3$
$$
    \psi_\la(n)
:=  \frac{n}{n+2} \frac{(n+2\la)(n+2\la+8)}{(n+2\la+1)^2}
\ge \frac{n}{n+2} \frac{n(n+8)}{(n+1)^2}
\ge 1\,.
$$
This proves the lower estimates \rf[B>].
\qed


\bigskip
For the proof of Theorem \ref{thm5}, we need yet one more lower
bound.

\begin{lemma}
For all $n \in \N$ and $\la > -\frac{1}{2}$, we have
\be \lb{c1>}
    [c_n(\la)]^2
\ge \frac{2}{2\la+1} n(n+\la)(n+2\la)\,.
\ee
\end{lemma}

\proof For any symmetric matrix $\mathbf{C} \in \R^{m\times m}$, its
largest eigenvalue $\mu_{\max}(\mathbf{C})$ satisfies the inequality
$\mu_{\max}(\mathbf{C})=\sup_{\Vert
\mathbf{x}\Vert=1}(\mathbf{C}\mathbf{x},\mathbf{x}) \ge
(\mathbf{C}\mathbf{e}_i,\mathbf{e}_i)= c_{ii}$, $1\leq i\leq m$.
Therefore,
$$
   [c_n(\la)]^2
 = \mu_{\max}(\B_n) \ge b_{mm}
 = 4 a_{mm}
$$
and by \rf[al]-\rf[wtal], with $f_0(x) = x(x+\frac{\la}{2})(x+\la)$, we have
$$
   4a_{mm} 
 = \frac{16}{2\la+1} f_0\Big(\frac{n}{2}\Big)
 = \frac{2}{2\la+1} n(n+\la)(n+2\la)\,.
$$
\qed

We will prove Theorem \ref{thm5} by establishing a slightly stronger
statement.

\begin{theorem}
For $n \ge 3$ and $\la > 2$, we have
\be \lb{e4'}
    \frac{1}{8} F(n,\la)
\le [c_n(\la)]^2
\le F(n,\la)
\ee
where
\be \lb{F1}
     F(n,\la) = \frac{n(n+\la)(n+2\la)(n+3\la)}{(\la+1)(\la+2)}
\ee
\end{theorem}

\proof 1) For the upper bound, using the upper bound in \rf[B1'], we
have \baa
      [c_n(\la)]^2
&\le& \frac{ n(n+\la)(n+2\la)(n+3\la+6)}{(\la+2)(\la+3)}
 =:  F(n, \la) \phi(n,\la)
\eaa
where
$$
   \phi(n,\la)
:= \frac{\la+1}{\la+3}\cdot \frac{n+3\la+6}{n+3\la} \le
\frac{\la+1}{\la+3} \cdot\frac{3 + 3\la + 6}{3+3\la}
 = 1\,, \qquad n \ge 3.
$$

2) For the lower bound, we consider two cases.

2a)~~ If $n > 5\la$, we use the lower estimate \rf[e1]
$$
    [c_n(\la)]^2
\ge \frac{1}{4} \frac{n^2(n+\la)^2}{(\la+1)(\la+2)}
=: \frac{1}{4} F(n,\la) \psi_1(n,\la) \,,
$$
where
$$
   \psi_1(n,\la)
:= \frac{n(n+\la)}{(n+2\la)(n+3\la)}
 = \frac{1}{(1 + \frac{2\la}{n})(1 + \frac{2\la}{n+\la})}
 > \frac{1}{(1 + \frac{2}{5})(1 + \frac{2}{6})}
 = \frac{5}{7}\cdot\frac{6}{8}
 > \frac{1}{2}
$$

2b) For $n \le 5\la$, we use the estimate \rf[c1>],
$$
    [c_n(\la)]^2
\ge \frac{2}{2\la+1} n(n+\la)(n+2\la)
\ge \frac{1}{\la+1} n(n+\la)(n+2\la)
 =  F(n,\la) \psi_2(n,\la)\,.
$$
where
$$
   \psi_2(n,\la)
:= \frac{\la+2}{n+3\la}
 > \frac{\la}{n+3\la}
\ge \frac{1}{5+3} = \frac{1}{8}\,.
$$
\qed

\noindent {\bf Proof of Theorem \ref{thm5}.} Since
$$
   \frac{3}{4} (n+2\la)^2 < (n+\la)(n+3\la) < (n+2\la)^2
$$
and
$$
    \frac{1}{\frac{3}{2}\la^2}
 <  \frac{1}{(\la+1)(\la+2)}
 <  \frac{1}{\la^2}, \qquad \la \ge 7,
$$
we derive from \rf[e4'] that
$$
   \frac{1}{16}\,\frac{n(n+2\la)^3}{\la^2}
 < [c_n(\la)]^2 < \frac{n(n+2\la)^3}{\la^2}\,,\quad \lambda\geq 7\,,
$$
and that proves \rf[e4].
\qed

\noindent {\bf Proof of Corollary \ref{thm6}.} Claim i) is
equivalent to
$$
    n \le \lim_{\la\to\infty} \frac{c_n(\la)^2}{2\la} \le 3n\,.
$$
The upper estimate follows from  \rf[e4'], while the lower estimate
follows from \rf[c1>]. Claim ii) follows from estimates \eqref{e2}.
\qed

\begin{remark} \rm
The approach proposed here is applicable for derivation of tight two
sided estimates for the best constant in the Markov $L_2$ inequality
with the Laguerre weight $\,w_{\alpha}(x)=x^{\alpha} e^{-x}$. The
results will appear in a forthcoming paper.
\end{remark}

\noindent {\bf Acknowledgement.} This research was performed during
a three week stay of the authors in the Oberwolfach Mathematical
Institute in April, 2016, within the Research in Pairs Program. The
authors thank the Institute for hospitality and the perfect research
conditions.


\newpage\noindent
{\sc Geno Nikolov} \smallskip\\
Department of Mathematics and Informatics\\
Universlty of Sofia \\
5 James Bourchier Blvd. \\
1164 Sofia \\
BULGARIA \\
{\it E-mail:} {\tt geno@fmi.uni-sofia.bg}

\bigskip\bigskip\noindent
{\sc Alexei Shadrin} \smallskip\\
Department of Applied Mathematics and Theoretical Physics (DAMTP) \\
Cambridge University \\
Wilberforce Road \\
Cambridge CB3 0WA \\
UNITED KINGDOM \\
{\it E-mail:} {\tt a.shadrin@damtp.cam.ac.uk}

\end{document}